\documentclass{amsart}

\usepackage{verbatim,amssymb,amsmath,amscd,latexsym,amsbsy,mathrsfs} 
\usepackage{graphicx}
\input{xy}
\xyoption{all}

\newtheorem{thm}{Theorem}[section]
\newtheorem{cor}[thm]{Corollary}

\DeclareMathAlphabet{\mathpzc}{OT1}{pzc}{m}{it}
 
\DeclareMathOperator{\Hom}{Hom}

\newcommand{\QQ}{\mathbb{Q}}
\newcommand{\ZZ}{\mathbb{Z}}

\newcommand {\C} {{\mathbb C}}

\newcommand {\Z} {{\mathbb Z}}
\newcommand {\Q} {{\mathbb Q}}

\newcommand {\dt} {{\bullet}}
\newcommand {\G} {{\mathbb G}}

\newcommand {\OO} {{\mathcal O}}
\newcommand {\A} {\mathbb{A}}

 \newtheorem{lemma}{Lemma}[section]

\begin{document}
\title{ Beilinson-Hodge cycles on semiabelian varieties} 
\author{
        Donu Arapura    
}
\address{Department of Mathematics\\
Purdue University\\
West Lafayette, IN 47907\\
U.S.A.}
\email{arapura@math.purdue.edu}
\author{
Manish Kumar
}
\address{Department of Mathematics\\
 Michigan State University\\
 East Lansing, MI-48824\\
 U.S.A.}
\email{mkumar@math.msu.edu}
 \maketitle

Given  a smooth not necessarily proper complex variety $U$, Beilinson
\cite{beilinson} conjectured that all Hodge cycles in $H^*(U,\Q)$
come from motivic cohomology, or more precisely that the so called regulator map
$$reg:CH^i(U,j)\otimes \Q\to Hom_{MHS}(\Q(-i), H^{2i-j}(U,\Q))$$
from Bloch's higher Chow group \cite{bloch} is surjective.
This is a very natural and appealing statement which
includes  the usual Hodge conjecture. Unfortunately, it has turned out
that it is not true in this generallity, c.f. \cite[9.11]{jannsen}, \cite{kl}. 
There is presumably a restricted range of $(i,j)$ for which this conjecture is
viable. For instance the line $j=0$, which corresponds  to 
the usual Hodge conjecture, should lie in this set. 
Work of  Asakura and Saito \cite{as} suggests that
 the conjecture should also  hold when $i=j$. Following these authors,
we refer to this special case as the Beilinson-Hodge conjecture.

Our goal here is to prove the Beilinson-Hodge   conjecture when
$U$ is either a semiabelian variety
or a product of smooth curves. The method is based on the study of
invariants under the Mumford-Tate group. 

\section{Reduction lemma}

We recall \cite{bloch, levine} that given a variety $U$, Bloch  has
defined a bigraded abelian group
$\bigoplus CH^i(U,j)$. 
The elements are represented by certain codimension $i$ algebraic cycles on
$U\times \mathbb{A}^j$.
There are  products
$$CH^i(U,j)\times CH^p(U,q)\to  CH^{i+p}(U,j+q)$$
when $U$ is smooth. 
 A cycle $Z\subset U\times \mathbb{A}^j$, representing an element
 of $ CH^i(U,j)$, has a fundamental class in
$$H^{2i}(U\times \A^j,U\times \partial \A^j)(i) \cong H^{2i-j}(U)(i)$$ 
where $\partial \A^j$ is a union of the hyperplanes corresponding to
the faces of $\A^j$ when viewed as an algebraic  simplex.
This extends to a homomorphism
$$reg:CH^i(U,j)\to Hom_{MHS}(\Z(-i), H^{2i-j}(U,\Z))$$
 This description was indicated in \cite{bloch}.
 Other explicit  constructions of this map can be found in
\cite{klm}, and \cite[\S 1]{as} for the subgroup of decomposable
cycles. From these formulas, it is clear that the map respects
products, and  the special case
$$reg:CH^1(U,1) =\OO(U)^*\to \Hom_{MHS}(\ZZ(-1),H^1(U,\ZZ))\subset H^1(U,\Z(1))$$
is just the composition of the inclusion $\OO(U)^*\subset \OO^{an}(U)^*$ with the
connecting map associated to the exponential sequence.

It is convenient to define the space of   Beilinson-Hodge cycles
$$BH^q(U) =  Hom_{MHS}(\Q(-q), H^{q}(U,\Q))$$
Then the  Beilinson-Hodge conjecture asserts that $CH^q(U,q)$ surjects onto
$BH^q(U)$. Note that the conjecture is only interesting for open
varieties, because it is vacuously true if the
variety is proper, since $BH^*=0$ in this case by
\cite{deligne-hodge}. 
The first nontrivial  case of the conjecture, when $q=1$,  turns out to be
 easy to understand and  prove, even integrally. It is not
 unreasonable to attribute this to Abel, since it is closely related
 to  his classical theorem.

\begin{thm}[Abel]\label{thm:abel}
  For any smooth variety $U$, the map
$$reg:\OO(U)^*\to \Hom_{MHS}(\ZZ(-1),H^1(U,\ZZ))$$
is surjective
\end{thm}

\begin{proof}
 Choose a smooth compactification $X$ such that $D=X-U$ has
 normal crossings. Let $d\OO^{an}_U$ denote the image of $d:\OO^{an}_U\to \Omega_U^{an1}$ in the category of sheaves.
 The group $H^1(U,\Z(1))$ is torsion free by the universal coefficient
 theorem, so it can be viewed as a subgroup of $H^1(U,\C)$.
An element in $H^1(U,\Z(1))$ is in $BH^1(U)$ if
  and only if  it lies in $F^1H^1(U)= \ker[H^1(U,\C)\to H^1(X,\OO_X)]$.
Chasing  the  following commutative diagram, with exact rows,
$$
\xymatrix{
H^0(U,\OO_U^{an*})\ar[r]^{\delta}\ar[d]^{d\log} & H^1(U,\Z(1))\ar[r]\ar[d] & H^1(U,\OO^{an}_U)\ar@{=}[d]\\
H^0(U,d\OO^{an}_U)\ar[r] & H^1(U,\C)\ar[r]       & H^1(U,\OO^{an}_U)\\
H^0(X,\Omega_X^1(\log D))\ar[r]\ar[u]& H^1(U,\C)\ar@{=}[u]\ar[r] & H^1(X,\OO_X)\ar[u]
}
$$
shows that the set of these classes coincides with $\{\delta(f)\mid d\log(f)\in H^0(\Omega_X(\log D))\}$.
The condition $d\log(f)\in H^0(\Omega_X(\log D))$ can be seen to force $f$ to have singularities of finite order
along $D$. Thus
 $$BH^1(U)\cap H^1(X,\Z)=\delta( \OO(U)^*).$$
\end{proof}

\begin{lemma}\label{lemma:key}
 If the products $BH^1(U)\times \ldots \times BH^1(U)\to BH^q(U)$  are surjective
 for all $q$, then the Beilinson-Hodge conjecture holds for $U$.
\end{lemma}

\begin{proof}
This follows from the following commutative diagram and theorem 1.1
$$
\xymatrix{
 CH^1(U,1)\times \ldots \times CH^1(U,1)\ar[r]\ar[d] & CH^q(U,q)\ar[d] \\ 
 BH^1(U)\times \ldots \times BH^1(U)\ar[r] & BH^q(U)
}
$$  
\end{proof}

\section{Mumford-Tate groups}

The category of rational mixed Hodge structures form a neutral
Tannakian category over $\Q$ \cite[chap II]{dm}.  Let $\langle H\rangle$ denote the
Tannakian category generated by a mixed Hodge structure $H$. This is the full subcategory
consisting of all subquotients
of tensor powers $T^{m,n}H=H^{\otimes m}\otimes (H^*)^{\otimes
  n}$. This construction extends to any set of Hodge structures.
The Mumford-Tate group $MT(H)$ is the group of tensor automorphisms
of the forgetful functor from $\langle H\rangle$ to $\Q$-vector
spaces.  By Tannaka duality $\langle H\rangle$ is equivalent to the
category of representations of this group.  When $H$ is a pure Hodge
structure, $MT(H)$ can be defined in a more elementary fashion as the
smallest $\Q$-algebraic group whose real points contains the image of
the torus defining the Hodge structure.  We define two auxillary
groups. The extended Mumford-Tate
group $EMT(H)$ is $MT(\langle H, \QQ(1)\rangle)$, and it surjects onto $MT(H)$. (Some authors
consider $EMT(H)$  to be the
Mumford-Tate group). The special Mumford-Tate group 
$SMT(H) = \ker[EMT(H)\to \G_m]$ with respect to the map that is induced by the inclusion 
$\langle \Q(1)\rangle\subset \langle H, \QQ(1)\rangle$.

\begin{thm}\label{thm:charofMT}
  \begin{enumerate}

  \item[]

  \item If $\Q(1)$ (respectively $\Q(m)$ with $m\not=0$) lies in
    $\langle H\rangle$, then $MT(H)$ is isomorphic (respectively
    isogenous) to $EMT(H)$. Otherwise $EMT(H)\cong MT(H)\times \G_m$.

  \item $MT(H)\subset GL(H)$ is the largest subgroup leaving every
    rational element of type $(0,0)$ in $T^{m,n}H$ invariant for all
    $m,n$. $SMT(H)$ leaves rational elements of type $(q,q)$ in
    $T^{m,n}H$ invariant for all $m,n,q$.

  \item If $H$ is pure and polarizable, then $MT(H)$ is connected and
    reductive.

  \item Let $H^{split} = \bigoplus_k Gr_k^WH$, then $MT(H)$ is a
    semidirect product of $MT(H^{split})$ with a unipotent group.
  \end{enumerate}
\end{thm}

\begin{proof}
  For the first statement, see \cite[pp 466-467]{milne}.  The next two
  properties are standard and proved in \cite[chap
  I]{dm}, although \cite[\S 2]{andre}  would be a more concise reference. The last part is
  essentially given in \cite{andre}. We indicate the proof for completeness.  Let $P$
  be the group linear automorphisms of $H$ preserving the flag
  $W_\dt$. The unipotent radical $UP\subset P$ is the subgroup which
  acts trivially on $Gr^W_k$.  We have inclusion of tensor categories
$$\iota:\langle H^{split} \rangle \to \langle H \rangle$$
with a right inverse $H' \mapsto (H')^{split}$. Therefore we get a split
surjection of Tannaka duals $\iota^*:MT(H)\to MT(H^{split})$. The
kernel $\iota^*$ lies in $UP$, and is therefore unipotent.
\end{proof}

\begin{cor}
$MT(H^{split})$ is the quotient of $MT(H)$ by its unipotent radical.  
\end{cor}

Let us turn to the case where $U$ is either a semiabelian variety
or a smooth curve. 
Set $MT(U) = MT(H^1(U))=EMT(H^1(U))$, where the last equality follows
from the theorem. Also let $SMT(U) =
SMT(H^1(U))$. 

Let $H= H^1(U)$ and let $W = W_1H=H^1(X)$. Choose a complementary subspace
$V$ to $W$ in $H$. We also know that $MT(U)$ preserves the weight filtration
on $H^1(U)$ (\cite[Lemma 2c]{andre}). Hence $\Phi$ the kernel of $MT(H)\to
MT(H^{split})$ and the unipotent radical of $MT(U)$ is a subspace of
$\Hom_{\Q}(V, W)$.

  \begin{cor}\label{cor:SMT}
As a subgroup of $GL(H) = GL(V\oplus W)$
    $$ SMT(U)=
\{\begin{pmatrix}
 I & 0 \\
 f & S 
 \end{pmatrix} \mid S\in SMT(W) \text{ and $f\in \Phi $}\}.$$
  \end{cor}

\section{Main theorem}

Let $H$ be the first cohomology of a semiabelian
variety or a smooth affine curve. We want to refine the description 
of $SMT(H)$ given by corollary \ref{cor:SMT}.
 We define three subspaces $V_i\subset
H$. Let
$V_3 = W_1H$, let $V_1 \subseteq H^{SMT(H)}$ be a
complement to $V_3$ in $ W_1H+H^{SMT(H)} $, and finally choose $V_2$ to be
a complement to $V_1+V_3$ in $H$. Thus we have a decomposition
\begin{equation}
  \label{eq:decompH}
  H = V_1\oplus V_2\oplus V_3  
\end{equation}
 with respect to which 
$SMT(H)$ becomes a subgroup of the following matrix group:
$$
\{\begin{pmatrix}
I& 0 & 0\\
 0 &I& 0\\
 0 & f & S 
 \end{pmatrix} \mid S\in SMT(V_3) \text{ and $f\in Hom(V_2,V_3)$}\}.$$
The unipotent radical $U(SMT(H))$  lies in the subgroup
\begin{equation}
  \label{eq:USMT}
  \{\begin{pmatrix}
I& 0 & 0\\
 0 &I& 0\\
 0 & f & I 
 \end{pmatrix} \mid S\in SMT(V_3) \text{ and $f\in Hom(V_2,V_3)$}\}.
\end{equation}

 \begin{lemma}\label{lemma:Ureduction}
   For any nonzero $u\in V_2$, we can
find a  $g\in U(SMT(H)) $ such that $gu\not=u$, or equivalently such that
$f(u)\not=0$ with respect to the matrix (\ref{eq:USMT}).
 \end{lemma}

 \begin{proof}
   Given a nonzero  $u\in V_2$, we have $g_1u\not= u$ for some $g_1\in
   SMT(H)$. Writing 
$$
g_1 = \begin{pmatrix}
I& 0 & 0\\
 0 &I& 0\\
 0 & f & S 
 \end{pmatrix}
$$
we see that $f(u)\not=0$. Set
$$
g_2=
\begin{pmatrix}
I& 0 & 0\\
 0 &I& 0\\
 0 & 0 & S^{-1} 
 \end{pmatrix}
$$
This lies in $SMT(H)$, since the map $SMT(H)\to SMT(H)/U(SMT(H))$
splits. Then $g= g_1g_2$ has the desired property.
 \end{proof}

Let 
$$BH^q(H) = \Hom(\Q(-q), H^{\otimes q})$$
for $H$ as above.

 \begin{thm}
   The product maps
    $BH^1(H)\times \ldots \times BH^1(H)\to BH^q(H)$  
    are surjective for all $q$
 \end{thm}

 \begin{proof}
To simplify book keeping,
we will usually write  tuples $(j_1,\ldots j_n)$ as strings $j_1\ldots j_n$.
Juxtaposition is used to denote concatenation of
strings, with exponents used for repetition. For example, $1^2\,2\,3^0= 1\,1\,2$.

(\ref{eq:decompH}) leads to a decomposition 
\begin{equation}
  \label{eq:decompHn}
  H^{\otimes n}=\bigoplus_{j_1,\ldots, j_n} V({j_1\ldots j_n}),
\end{equation}
where
$$ V(j_1\ldots j_n)=V_{j_1}\otimes\ldots\otimes V_{j_n}$$

   Let  $\tau \in BH^n(H)$ i.e. suppose that it is a
Beilinson-Hodge cycle. Our goal is
to show that $\tau\in BH^1(H)^{\otimes n}$. 
Let us decompose
$$\tau =\sum \tau_{j_1\ldots j_n}$$
with respect to \eqref{eq:decompHn}.
It suffices to show that $\tau \in V_1^{\otimes
  n} $, since $V_1\subseteq BH^1(H)$. After replacing $\tau$ by $\tau-\tau_{1^n}$,
we will show $\tau$ equals $0$.

We next argue that any component $\tau'=\tau_{j_1j_2\ldots j_n}$
with all of the $j_i\in \{1,2\}$  must be zero. 
Assume that $\tau'\not=0$, then we will derive a contradiction.
Let
$$\tau_{j_1j_2\ldots j_n}=
x_1\otimes x_2\otimes \ldots\otimes x_n$$
with $x_i\in V_{j_i}$.
From the previous paragraph, $j_1\ldots j_n = 1^{n_1}2^{n_2}\,1^{n_3}\ldots$ must have
at least one $2$. Since $u = x_{n_1+1} \in V_2-\{0\}$, we can 
choose a  $g\in U(SMT(H))$ so that $f(u)\not= 0$, with $f$ as in \eqref{eq:USMT}. Then
$g\tau' -\tau'$ will have a nonzero component in $V({1^{n_1\,}3\,
  2^{n_2-1} 1^{n_3}\ldots })$.
We must have $g\tau-\tau=0$, since $\tau$ is invariant under $SMT(H)$
by theorem~\ref{thm:charofMT}.
Thus $\tau$ must have another term $\tau''$ whose  image under $g-I$ has
 a nonzero component of   type $1^{n_1}\,3\,2^{n_2-1}\ldots$.
The only possible candidate is  $\tau''=\tau_{1^{n_1}32^{n_2-1}\ldots}$.
However, after
expanding this as a product of $x_i$'s as above, we can see that
$(g-I)\tau''$ has no nonzero components of the required type.
For example,
$(g-I)\tau''=0$ if the second $2$ is absent from $j_1j_2\ldots j_n$,
$(g-I)\tau''$ is sum of types $1^{n_1}3^21^{n_3}, 1^{n_1}321^{n_3}$ and $1^{n_1}231^{n_3}$ if $j_1j_2\ldots
j_n=1^{n_1}2^21^{n_3}$ and so on. Therefore $\tau'=0$ as claimed.

To conclude, we note that
the projection of a nonzero Beilinson-Hodge cycle to $(Gr^W_2 H)^{\otimes n}$
must be nonzero.
We deduce from the previous paragraph that for every component of
$\tau$, must have at least one $j_i=3$. This implies that  $\tau$ 
projects to zero in $(Gr^W_2 H)^{\otimes n}$. Therefore it must
already  be zero.
 \end{proof}

\begin{cor}
  The Beilinson-Hodge conjecture holds for a product of smooth curves.
\end{cor}

\begin{proof}
  Let $U = \prod U_i$, where $U_i$ are smooth curves. Let $H= H^1(U)$.
  Then by K\"unneth's
  formula and the theorem, the conditions of lemma~\ref{lemma:key} hold.
\end{proof}

\begin{cor}
   The Beilinson-Hodge conjecture holds for a semiabelian variety.
\end{cor}

\begin{proof}
 Let $U$ be a semiabelian variety. Let $H = H^1(U)$. By
the theorem, we have that $BH^n(H) = BH^1(H)^{\otimes
  n}$. Now observe that $H^*(U) = \wedge^* H$ which is a direct
summand of the tensor algebra. So the BH cycles on $H^n(U)$ are given
by products of BH-cycles on $H$.
\end{proof}

The referee pointed out the following interesting corollary  which can be proved along the
same lines as the first corollary.

\begin{cor}

  Let $U=\prod U_i$  be a product of $n$ smooth  curves with smooth projective completions $X_i$.
Then $BH^n(U) \not=0$ if and only if there exists  torsion cycles in $J(X_i)$ with nonempty support
on $X_i-U_i$ for each $i$.
\end{cor}

\begin{proof}
Using the theorem, this can be reduced to the case of $n=1$.
By theorem~\ref{thm:abel}, a  nonzero element of $BH^1(U_1)$ lifts
to an element  $f\in \OO(U_1)^*\otimes \Q$, which in turn defines a divisor $(f)\in Div(U_1)\otimes \Q$ 
with nonempty support in $X_1-U_1$. Conversely, any such $\Q$-divisor determines a
nonzero element of $BH^1(U_1)$
\end{proof}

\section*{Acknowledments}

We would like to thank
Jinhyun Park for some useful conversations about higher Chow groups,
and the referee for some helpful remarks.

The first author was partially supported by the NSF.

\end{document}